\documentclass[11pt]{article}
\usepackage{tikz}
\usepackage[all]{xy}
\usepackage{amsfonts}
\usepackage{amsfonts}
\usepackage{amsfonts}
\usepackage{mathrsfs}
\usepackage{bm}
\usepackage{cite}
\usepackage[colorlinks,linkcolor=blue,citecolor=blue]{hyperref}
\usepackage{amssymb,amsmath} 
\usepackage{makecell}
\usepackage{tikz}
 \allowdisplaybreaks
\usetikzlibrary{arrows,shapes,chains}
 \allowdisplaybreaks

\catcode`!=11
\let\!int\int \def\int{\displaystyle\!int}
\let\!lim\lim \def\lim{\displaystyle\!lim}
\let\!sum\sum \def\sum{\displaystyle\!sum}
\let\!sup\sup \def\sup{\displaystyle\!sup}
\let\!inf\inf \def\inf{\displaystyle\!inf}
\let\!cap\cap \def\cap{\displaystyle\!cap}
\let\!max\max \def\max{\displaystyle\!max}
\let\!min\min \def\min{\displaystyle\!min}
\let\!frac\frac \def\frac{\displaystyle\!frac}
\catcode`!=12

\let\oldsection\section
\renewcommand\section{\setcounter{equation}{0}\oldsection}

\allowdisplaybreaks
\def\pf{\it{Proof.}\rm\quad}

\def\N{\mathbb{N}}\def\Z{\mathbb{Z}}

\def\sn{\sum\limits_{k=1}^n}

\def\su{\sum\limits_{n=1}^\infty}
\def\sk{\sum\limits_{k=1}^\infty}
\def\sj{\sum\limits_{j=1}^\infty}
\def\t{\widetilde{t}}

\def\z{\zeta}

\def\a{^{(A)}}
\def\B{^{(B)}}
\def\C{^{(C)}}
\def\ab{^{(AB)}}
\def\abc{^{(ABC)}}

\newtheorem{defn}{Definition}[section]
\newtheorem{thm}{Theorem}[section]
\newtheorem{lem}[thm]{Lemma}
\newtheorem{cor}[thm]{Corollary}
\newtheorem{con}[thm]{Conjecture}

\newtheorem{exa}{Example}[section]

\begin{document}
\title {\bf Explicit Evaluations for Several Variants of Euler Sums}
\author{
{Ce Xu\thanks{Email: 19020170155420@stu.xmu.edu.cn}}\\[1mm]
\small{$*$ School of Mathematics and Statistics, Anhui Normal University,}\\ \small{Wuhu 241000, P.R. China}}

\date{}
\maketitle \noindent{\bf Abstract} We study several variants of Euler sums by using the methods of contour integration and residue theorem. These variants exhibit nice properties such as closed forms, reduction, etc., like classical Euler sums. In addition, we also define a variant of multiple zeta values of level 2, and give some identities on relations between these variants of Euler sums and the variant of multiple zeta values.
\\[2mm]


\section{Introduction}

The classical Euler sums was introduced and studied by Flajolet and Salvy \cite{FS1998}, which are defined by the infinite series
\begin{align}\label{a1}
{S_{{\bf p},q}} := \sum\limits_{n = 1}^\infty  {\frac{{H_n^{\left( {{p_1}} \right)}H_n^{\left( {{p_2}} \right)} \cdots H_n^{\left( {{p_r}} \right)}}}
{{{n^q}}}},
\end{align}
where $H_n^{(p)}$ stands for the $p$-th generalized harmonic number, which is defined by
\begin{align*}
H_n^{(p)}: = \sum\limits_{k=1}^n {\frac{1}{{{k^p}}}},\quad H_n\equiv H_n^{(1)}\quad {\rm and}\quad H_0^{(p)}:=0.
\end{align*}
In this paper we consider the related quantities
\begin{align}\label{a2}
{R_{{\bf p},q}} := \sum\limits_{n = 0}^\infty  {\frac{{H_n^{\left( {{p_1}} \right)}H_n^{\left( {{p_2}} \right)} \cdots H_n^{\left( {{p_r}} \right)}}}
{{{(n+1/2)^q}}}},
\end{align}
which we call Euler $R$-sums. Here ${\bf p}:=(p_1,p_2,\ldots,p_r)\ (r,p_i\in \N,\ i=1,2,\ldots,r)$ with $p_1\leq p_2\leq \cdots\leq p_r$ and $q\geq 2$. The quantity $w:={p _1} +  \cdots  + {p _r} + q$ is called the weight and the quantity $r$ is called the degree (order). For convenience, repeated summands in partitions are indicated by powers, so that for instance
\[{S_{{1^2}{2^3}4,q}} = {S_{112224,q}} = \sum\limits_{n = 1}^\infty  {\frac{{H_n^2[H^{(2)} _n]^3{H^{(4)} _n}}}{{{n^q}}}}\quad\text{and} \quad R_{1^32^25,q}=R_{111225,q}=\sum_{n=0}^\infty \frac{H_n^3[H_n^{(2)}]^2H_n^{(5)}}{(n+1/2)^q}. \]

In particular, in (\ref{a1}), we put a bar on top of $q$ if there is a sign $(-1)^{n-1}$ appearing in the
denominator on the right; In (\ref{a2}), we put a bar on top of $q$ if there is a sign $(-1)^{n}$ appearing in the
denominator on the right, namely,
\begin{align}\label{ab1}
{S_{{\bf p},{\bar q}}} := \sum\limits_{n = 1}^\infty  {\frac{{H_n^{\left( {{p_1}} \right)}H_n^{\left( {{p_2}} \right)} \cdots H_n^{\left( {{p_r}} \right)}}}
{{{n^q}}}}(-1)^{n-1}
\end{align}
and
\begin{align}\label{ab2}
{R_{{\bf p},{\bar q}}} := \sum\limits_{n = 0}^\infty  {\frac{{H_n^{\left( {{p_1}} \right)}H_n^{\left( {{p_2}} \right)} \cdots H_n^{\left( {{p_r}} \right)}}}
{{{(n+1/2)^q}}}}(-1)^n.
\end{align}

Similarly, we define the alternating harmonic number ${\bar H}^{(p)}_n$ by
\begin{align*}
{ \bar H}_n^{(p)}: = \sum\limits_{k=1}^n {\frac{(-1)^{k-1}}{{{k^p}}}},\quad {\bar H}_n\equiv {\bar H}^{(1)}_n \quad {\rm and}\quad {\bar H}_0^{(p)}:=0,
\end{align*}
which also was introduced in \cite{FS1998}. In (\ref{a1})-(\ref{ab2}), if replace ``$H^{(p_j)}_n$" by ``${\bar H}^{(p_j)}_n$" in the numerator of the summand, we put a ``bar'' on the top of $p_j$. For example,
\begin{align*}
S_{{\bar p}_1{p}_2p_3{{\bar p}_4},{q}}:=\su \frac{{\bar H}^{(p_1)}_n{H}^{(p_2)}_nH^{(p_3)}_n{\bar H}^{(p_4)}_n}{n^q},\quad S_{p_1{\bar p}_2p_3{{\bar p}_4},{\bar q}}:=\su \frac{H^{(p_1)}_n{\bar H}^{(p_2)}_nH^{(p_3)}_n{\bar H}^{(p_4)}_n}{n^q}(-1)^{n-1}
\end{align*}
and \begin{align*}
R_{{\bar p}_1{p}_2{\bar p}_3{{\bar p}_4},{q}}:=\sum_{n=0}^\infty \frac{{\bar H}^{(p_1)}_n{H}^{(p_2)}_n{\bar H}^{(p_3)}_n{\bar H}^{(p_4)}_n}{(n+1/2)^q},\quad R_{{\bar p}_1{\bar p}_2p_3{{\bar p}_4},{\bar q}}:=\sum_{n=0}^\infty \frac{{\bar H}^{(p_1)}_n{\bar H}^{(p_2)}_nH^{(p_3)}_n{\bar H}^{(p_4)}_n}{(n+1/2)^q}(-1)^{n}.
\end{align*}
The sums of types above (one of more the $p_j$ or $q$ barred) are called the alternating Euler sums and alternating Euler $R$-sums, respectively.  Classical Euler sums may be studied through a profusion of methods: combinatorial, analytic and algebraic. There are many other researches on Euler sums and Euler type sums. Some related results for Euler sums may be seen in the works of \cite{BBG1994,BBG1995,F2005,M2014,W2017,XW2019,Z2019} and references therein.
Clearly, (alternating) Euler sums can be expressed in terms of a $\mathbb{Q}$-linear combination of (alternating) multiple zeta values (explicit formula of (alternating) Euler sums via (alternating) multiple zeta values can be found in \cite{XW2018}).
Flajolet and Salvy \cite{FS1998} also pointed out that every Euler sum is a $\mathbb{Q}$-linear combination of multiple zeta values, but not explicit formula. For non-zero integers $k_1,k_2,\ldots,k_r$ and $k_1\neq 1$, the (alternating) multiple zeta values (abbr. MZVs) are defined by
\begin{align}\label{a3}
&\zeta( \mathbf{k})\equiv\zeta(k_1, \ldots, k_r):=\sum\limits_{n_1>\cdots>n_r\geq 1}\prod\limits_{j=1}^r n_j^{-|k_j|}{\rm sgn}(k_j)^{n_j},
\end{align}
where for convergence $|k_1|+\cdots+|k_j|> j$ for $j= 1, 2, \ldots, r$, and
\[{\rm sgn}(k_j):=\begin{cases}
   1  & \text{\;if\;} k_j>0,  \\
   -1, & \text{\;if\;} k_j<0.
\end{cases}\]
Here, we call $l(\mathbf{k}):=r$ and $w:=|k_1|+|k_2|+\cdots+|k_r|$ the depth and the weight of (alternating) multiple zeta values, respectively. The study of multiple zeta values began in the early 1990s with the works of Hoffman \cite{H1992} and Zagier \cite{DZ1994}. The study of multiple zeta values have attracted a lot of
research in the area in the last three decades.
For detailed history and applications, please see \cite{BBBL1997,BBV2010,BBBL2001} and the book of Zhao \cite{Z2016}. Recently, several variants of multiple zeta values were introduced and studied by Hoffman, Kaneko and Tsumura. Let $k_1>1,k_2,\ldots,k_r$ be positive integers. Recently,
Hoffman \cite{H2019} introduced and studied the multiple $t$-values (abbr. MtVs)
\begin{align}
t(k_1,k_2,\ldots,k_r):&=\sum\limits_{n_1>\cdots>n_r\geq 1\atop n_i \ {\rm odd}} \frac{1}{n_1^{k_1}n_2^{k_2}\cdots n_k^{k_r}}\nonumber\\&=\sum\limits_{n_1>\cdots>n_r\geq 1} \frac{1}{(2n_1-1)^{k_1}(2n_2-1)^{k_2}\cdots (2n_r-1)^{k_r}}.
\end{align}
As a contrast, Kaneko and Tsumura \cite{KTA2019} introduced and studied the multiple $T$-values (abbr. MTVs)
\begin{align}
T(k_1,k_2,\ldots,k_r):&=2^r \sum_{n_1>n_2>\cdots>n_r>0\atop n_i\equiv r-i+1\ {\rm mod}\ 2} \frac{1}{n_1^{k_1}n_2^{k_2}\cdots n_r^{k_r}}\nonumber\\
&=2^r\sum\limits_{n_1>n_2>\cdots>n_k>0} \frac{1}{(2n_1-r)^{k_1}(2n_2-r+1)^{k_2}\cdots (2n_r-1)^{k_r}}.
\end{align}
As their normalized versions,
\begin{align*}
&\widetilde{t}(k_1,k_2,\ldots,k_r):=2^{k_1+\cdots+k_r}t(k_1,k_2,\ldots,k_r),\\
&\widetilde{T}(k_1,k_2,\ldots,k_r):=2^{k_1+\cdots+k_r-r}T(k_1,k_2,\ldots,k_r).
\end{align*}
It is clear that these values can be written as a linear combination of alternating multiple zeta values.

In this paper, we define the following variant of the multiple zeta value of level 2,
\begin{align}\label{a4}
R(k_1,k_2,\ldots,k_r):=2^{|k_1|+|k_2|+\cdots+|k_r|} \sum_{n_1>\cdots>n_r>0} \frac{{\rm sgn}(k_1)^{n_1}{\rm sgn}(k_2)^{n_2}\cdots {\rm sgn}(k_r)^{n_r}}{(2n_1-1)^{|k_1|} (2n_{2})^{|k_2|}\cdots(2n_r)^{|k_r|}},
\end{align}
for non-zero integers $k_1,k_2,\ldots,k_r$ with $k_1\neq 1$, which we call (alternating) multiple $R$-values (abbr. MRVs), if $r=1$ and $k_1=k\geq 2$, we let
$$R(k):=2^{k}\sum_{n=1}^\infty \frac{1}{(2n-1)^k}\quad (k\geq 2).$$
Hence, $R(k)=\t(k)=\widetilde{T}(k)=(2^k-1)\z(k)$ for $k\geq 2$. For convenience, we let $R(1):=2\log(2)$.

In (\ref{a3}) and (\ref{a4}), we may compactly indicate the presence of an alternating sign. When ${\rm sgn}(k_j)=-1$,  by placing a bar over the
corresponding integer exponent $|k_j|$. For example, we write
\begin{align*}
{\zeta}( {\bar 2,3,\bar 1,4} )={\zeta}( {-2,3,- 1,4})\quad \text{and}\quad R( {\bar 3,2,\bar 4,1})={R}( {-3,2,-4,1}).
\end{align*}

Obviously, from definition, non-alternating MRVs (namely $\forall k_j\geq 1$ and $k_1>1$) can be written as a linear combination of alternating multiple zeta values. For all $k_j\geq 1$ with $k_1>1$, it is easy to deduce that
\begin{align*}
R(k_1,k_2,\ldots,k_r)=2^{k_1+\cdots+k_r-r} \sum_{\sigma_j\in \{\pm 1\}\atop j=2,3,\ldots,r} \left\{\z(k_1,\sigma_2k_2,\ldots,\sigma_rk_r)-\z({\bar k_1},\sigma_2k_2,\ldots,\sigma_rk_r) \right\}.
\end{align*}
According to the definitions of Euler $R$-sum and multiple $R$-value, it is easy to see that every (alternating) Euler $R$-sum of weight $w$ and degree $r$ is a $\mathbb{Q}$-linear
combination of (alternating) MRVs of weight $w$ and depth at most $r+1$.

The purposes of this paper are to establish explicit evaluations of (alternating) Euler $R$-sums and MRVs for lower depths by using the methods of contour integration and residue theorem, which was developed in our previous paper \cite{X2020}.

The remainder of this paper is organized as follows.

In Section \ref{sec2}, we define some generic cotangent, digamma functions and basic notations and state some results about them given in our previous paper \cite{X2020}. Moreover, we give a residue theorem.

In Section \ref{sec3}, we use the results of the previous Section \ref{sec2} and residue theorem to evaluate Euler $R$-sums. In particular, we establish explicit formulas of linear and quadratic Euler $R$-sums. Further, we prove the reducible for general Euler $R$-sums.

In Section \ref{sec4}, we prove the generating function of the `height one' MRVs, and use the method of residue theorem to evaluate the MRVs of depth three.

\section{Definitions and Lemmas}\label{sec2}
We give three definitions, which were defined in \cite{X2020}, and give four lemmas. Let $A:=\{a_k\}\ (-\infty < k < \infty)$ be a sequence of complex numbers with ${a_k} = o\left( {{k^\alpha }} \right)\ (\alpha  < 1)$ if $k\rightarrow \pm \infty$. For convenience, we let $A_1$ and $A_2$ to denote the constant sequence $\{a_k\}=\{(1)^k\}$ and alternating sequence $\{a_k\}=\{(-1)^k\}$, respectively.

\begin{defn}\label{def1} With $A$ defined above, we define the parametric digamma function $\Psi \left( { - s;A} \right)$ by
\begin{align}\label{1.1}
\Psi \left( { - s;A} \right):= \frac{{{a_0}}}{s} + \sum\limits_{k = 1}^\infty  {\left( {\frac{{{a_k}}}{k} - \frac{{{a_k}}}{{k - s}}} \right)}\quad (s\in \mathbb{C}\setminus (\N\cup\{0\})).
\end{align}
\end{defn}
Clearly, if $A=A_1$, then the parametric digamma function $\Psi \left( { - s;A} \right)$ becomes the classical digamma function $\psi \left( { - s} \right)+\gamma$.

\begin{defn}\label{def2} Define the cotangent function with sequence A by
\begin{align}\label{1.2}
 \pi \cot \left( {\pi s;A} \right) &=  - \frac{{{a_0}}}{s} + \Psi \left( { - s;A} \right) - \Psi \left( {s;A} \right)\nonumber \\
  &= \frac{{{a_0}}}{s} - 2s\sum\limits_{k = 1}^\infty  {\frac{{{a_k}}}{{{k^2} - {s^2}}}}.
\end{align}
\end{defn}
It is clear that if letting $A=A_1$ or $A_2$ in (\ref{1.2}), respectively, then it become
\begin{align*}
&\cot \left( {\pi s;A_1} \right) = \cot \left( {\pi s} \right)\quad\text{and}\quad \cot \left( {\pi s;A_2} \right) = \csc \left( {\pi s} \right).
\end{align*}

We now provide notations that will be used throughout this paper.

\begin{defn} For nonnegative integers $j\geq 1$ and $n$, we define
\begin{align*}
&R\a(j):=\sk \frac{a_k}{(k-1/2)^j},\quad R\a(1):=\sk \left(\frac{a_k}{k-1/2}-\frac{a_k}{k}\right),\\
&\widehat{R}\a(j):=\sk \frac{a_{k-1}}{(k-1/2)^j},\quad \widehat{R}\a(1):=\sk \left(\frac{a_{k-1}}{k-1/2}-\frac{a_k}{k}\right),\\
&E\a_n(j):=\sn \frac{a_{n-k}}{k^j},\quad E\a_0(j):=0,\quad {\bar E}\a_n(j):=\sn \frac{a_{k-n-1}}{k^j},\quad {\bar E}\a_0(j):=0,\\
&F\a_n(j)= \left\{ {\begin{array}{*{20}{c}} \sk \frac{a_{k+n}-a_k}{k}
   {,\ j=1,}  \\
   {\sk \frac{a_{k+n}}{k^j},\ \ \ \ \ \;\;\;j>1,}  \\
\end{array} } \right.,\quad{\bar F}\a_n(j)= \left\{ {\begin{array}{*{20}{c}} \sk \frac{a_{k-n}-a_k}{k}
   {,\ j=1,}  \\
   {\sk \frac{a_{k-n}}{k^j},\ \ \ \ \ \;\;\;j>1,}  \\
\end{array} } \right.\\
&G\a_n(j):=E\a_n(j)-{\bar E}\a_{n-1}(j)-\frac{a_0}{n^j},\quad G\a_0(j):=0,\\
&L\a_n(j):=F\a_n(j)+(-1)^j{\bar F}\a_n(j),\quad M\a_n(j):=E\a_n(j)+(-1)^j F\a_n(j),\\
&{\bar M}\a_n(j):={\bar F}\a_n(j)-{\bar E}\a_{n-1}(j),\quad T\a_n(j):=G\a_n(j)+(-1)^j L\a_n(j).
\end{align*}
\end{defn}

It is clear that $M\a_n(j)+{\bar M}\a_n(j)=T\a_n(j)+\frac{a_0}{n^j}$, and
if $A=A_1$ and $A_2$, then
\begin{align*}
&M^{(A_1)}_n(j)=H^{(j)}_n+(-1)^j\z(j),\ {\bar M}^{(A_1)}_n(j)=\z(j)-H^{(j)}_{n-1},\ T^{(A_1)}_n(j)=(1+(-1)^j)\z(j),\\
&M^{(A_2)}_n(j)=(-1)^{n-1}{\bar H}^{(j)}_n+(-1)^j \left\{ {\begin{array}{*{20}{c}} (1-(-1)^n)\log(2)
   {,\ j=1,}  \\
   {(-1)^{n-1}{\bar \z}(j),\ \ \ \ \ \;\;\;j>1,}  \\
\end{array} } \right.\\
&{\bar M}^{(A_2)}_n(j)=(-1)^{n}{\bar H}^{(j)}_{n-1}+\left\{ {\begin{array}{*{20}{c}} (1-(-1)^n)\log(2)
   {,\ j=1,}  \\
   {(-1)^{n-1}{\bar \z}(j),\ \ \ \ \ \;\;\;j>1,}  \\
\end{array} } \right.\\
&T^{(A_2)}_n(j)=(-1)^{n-1}(1+(-1)^j){\bar \z}(j).
\end{align*}

Next, we give four lemmas, these lemmas are basic tools that will be used throughout this paper.
\begin{lem}(\cite{X2020})\label{lem2.1} Let $p\geq 1$ and $n$ be nonnegative integers, if $|s-n|<1$ with $s\neq n$, then
\begin{align}\label{b3}
\frac{{{\Psi ^{\left( {p - 1} \right)}}\left( { - s;A} \right)}}{{\left( {p - 1} \right)!}}=\frac{1}{{{{\left( {s - n} \right)}^p}}}\left\{a_n-\sj (-1)^j\binom{j+p-2}{p-1} M\a_n(j+p-1)(s-n)^{j+p-1} \right\}.
\end{align}
\end{lem}
\begin{lem}(\cite{X2020})\label{lem2.2} Let $p$ and $n$ be positive integers, if $|s+n|<1$, then
\begin{align}\label{b4}
\frac{{{\Psi ^{\left( {p - 1} \right)}}\left( { - s;A} \right)}}{{\left( {p - 1} \right)!}}=(-1)^p \sj \binom{j+p-2}{p-1} {\bar M}\a_n(j+p-1)(s+n)^{j-1}.
\end{align}
\end{lem}
\begin{lem}(\cite{X2020})\label{lem2.3} With $\cot(\pi s;A)$ defined above, if $|s-n|<1$ with $s\neq n\ (n\in \Z)$, then
\begin{align}\label{b5}
\pi \cot(\pi s;A)=\frac{a_{|n|}}{s-n}-\sj (-\sigma_n)^j T\a_{|n|}(j)(s-n)^{j-1},
\end{align}
where $\sigma_n$ is defined by the symbol of $n$, namely,
\begin{align*}
\sigma_n:= \left\{ {\begin{array}{*{20}{c}}\ 1
   {,\  n\geq 0,}  \\
   {-1,\ n<0.}  \\
\end{array} } \right.
\end{align*}
\end{lem}

We define a kernel function $\xi \left( s \right)$ by the two requirements: 1. $\xi \left( s \right)$ is meromorphic in the whole complex plane. 2. $\xi \left( s \right)$ satisfies $\xi \left( s \right)=o(s)$ over an infinite collection of circles $\left| s \right| = {\rho _k}$ with ${\rho _k} \to \infty $. Applying these two conditions of kernel function $\xi \left( s \right)$, Flajolet and Salvy gave the following residue lemma.
\begin{lem}(\cite{FS1998})\label{lem3.1}
Let $\xi \left( s \right)$ be a kernel function and let $r(s)$ be a rational function which is $O(s^{-2})$ at infinity. Then
\begin{align}\label{b6}
\sum\limits_{\alpha  \in O} {{\mathop{\rm Res}}{{\left[ {r\left( s \right)\xi \left( s \right)},s = \alpha  \right]}}}  + \sum\limits_{\beta  \in S}  {{\mathop{\rm Res}}{{\left[ {r\left( s \right)\xi \left( s \right)},s = \beta  \right]}}}  = 0.
\end{align}
where $S$ is the set of poles of $r(s)$ and $O$ is the set of poles of $\xi \left( s \right)$ that are not poles $r(s)$ . Here ${\mathop{\rm Re}\nolimits} s{\left[ {r\left( s \right)},s = \alpha \right]} $ denotes the residue of $r(s)$ at $s= \alpha$.
\end{lem}

\section{Explicit Evaluations of Euler $R$-sums}\label{sec3}
In this section, we study the reduction of (alternating) Euler $R$-sums by using contour integration. We need the following identities.
\begin{align}
&\lim_{s\rightarrow -1/2}\Psi^{(p-1)}(-s;A)=(-1)^p(p-1)!\widehat{R}\a(p),\label{c1}\\
&\lim_{s\rightarrow -1/2} \frac{d^m}{ds^m}(\pi \cot(\pi s;A))=m!\left((-1)^mR\a(m+1)-\widehat{R}\a(m+1)\right).\label{c2}
\end{align}
These identities can be proved by direct calculations. Let $B:=\{b_k\}\ (-\infty < k < \infty)$ be a sequence of complex numbers with ${b_k} = o\left( {{k^\beta }} \right)\ (\beta  < 1)$ if $k\rightarrow \pm \infty$.
\begin{thm}\label{thm3.1} For positive integers $p$ and $q>1$,
\begin{align}\label{c3}
&(-1)^{p+q}\su \frac{{\bar M}\B_n(p)}{(n-1/2)^q}a_n+\su \frac{M\B_{n-1}(p)}{(n-1/2)^q}a_{n-1}\nonumber\\
&=(-1)^p \sum_{j=1}^p \binom{p+q-j-1}{q-1} \su \frac{T\a_{n-1}(j)}{(n-1/2)^{p+q-j}}b_{n-1}\nonumber\\ &\quad -(-1)^p \sum_{k=0}^{q-1} \binom{p+q-k-2}{p-1} \left((-1)^k R\a(k+1)-\widehat{R}\a(k+1) \right)\widehat{R}\B(p+q-k-1)\nonumber\\
&\quad-(-1)^p \binom{p+q-1}{p} \widehat{R}\ab(p+q),
\end{align}
where $\widehat{R}\ab(q)$ is defined by
\[\widehat{R}\ab(q):=\sum_{k=1}^\infty \frac{a_{k-1}b_{k-1}}{(k-1/2)^q}.\]
\end{thm}
\pf In the context of this paper, the theorem results from applying the kernel function
\[\pi \cot \left( {\pi s};A \right)\frac{{{\Psi ^{\left( {p - 1} \right)}}\left( { - s;B} \right)}}{{\left( {p - 1} \right)!}}\]
to the base function $r(s)=(s+1/2)^{-q}$. Namely, we need to compute the residue of the function
\[{f_1}( {s;A,B}): = \pi \cot \left( {\pi s} ;A\right)\frac{{{\Psi ^{\left( {p - 1} \right)}}\left( { - s;B} \right)}}{{\left( {p - 1} \right)!{(s+1/2)^q}}}.\]
The only singularities are poles at the integers. At a negative integer $-n$ the pole is simple and the residue is
\[{\rm{Res}}\left[ {{f_1}( {s;A,B} ),s =  - n} \right] =(-1)^{p+q} \frac{{\bar M}\B_n(p)}{(n-1/2)^q}a_n.\]
At a non-negative integer $n$, the pole has order $p+1$ and the residue is
\begin{align*}
 {\rm Res}[f_1(s;A,B),s=n]&=(-1)^p \binom{p+q-1}{p} \frac{a_nb_n}{(n+1/2)^{p+q}} + \frac{{ M}\B_n(p)}{(n+1/2)^q}a_n\\&\quad-(-1)^pb_n \sum_{j=1}^p \binom{p+q-j-1}{q-1} \frac{T\a_n(j)}{(n+1/2)^{p+q-j}}.
\end{align*}
Finally the residue of the pole of order $q$ at $-1/2$ is found to be
\begin{align*}
&{\rm Res}[f_1(s;A,B),s=0]\\&=(-1)^p \sum_{k=0}^{q-1} \binom{p+q-k-2}{p-1} \left((-1)^k R\a(k+1)-\widehat{R}\a(k+1) \right)\widehat{R}\B(p+q-k-1).
\end{align*}
Here we used the formulas (\ref{b3})-(\ref{b5}) and (\ref{c1})-(\ref{c2}).
Summing these three contributions yields the statement of the theorem. \hfill$\square$

Clearly,
\begin{align*}
&R^{(A_1)}(j)=R(j),\quad R^{(A_1)}(1)=2\log(2)\quad {\rm and} \quad R^{(A_2)}(j)=R({\bar j}),\\
&\widehat{R}^{(A_1)}(j)=R(j),\quad \widehat{R}^{(A_1)}(1)=2\log(2)\quad {\rm and} \quad \widehat{R}^{(A_2)}(j)=-R({\bar j}).
\end{align*}
If setting $A,B\in\{A_1,A_2\}$ in Theorem \ref{thm3.1}, then we obtain the following results of four (alternating) linear Euler $R$-sums.
\begin{cor} For positive integers $p$ and $q\geq 2$,
\begin{align}
&\begin{aligned}\label{c4}
(1-(-1)^{p+q})R_{p,q}=&(-1)^p\sum_{j=0}^p (1+(-1)^j)\binom{p+q-j-1}{q-1}\zeta(j)R(p+q-j)\\
&+(-1)^p\sum_{k=0}^{q-1} (1-(-1)^k)\binom{p+q-k-2}{p-1}R(k+1)R(p+q-k-1)\\
&-(-1)^p(1+(-1)^q)\zeta(p)R(q),
\end{aligned}\\
&\begin{aligned}\label{c5}
(1-(-1)^{p+q})R_{\bar p,q}=&-(-1)^p\sum_{j=0}^p (1+(-1)^j)\binom{p+q-j-1}{q-1}\zeta(\bar j)R({p+q-j})\\
&-(-1)^p\sum_{k=0}^{q-1} (1+(-1)^k)\binom{p+q-k-2}{p-1}{R}(\overline{k+1}){R}(\overline{p+q-k-1})\\
&+(-1)^p(1+(-1)^q)\zeta(\bar p)R(q),
\end{aligned}\\
&\begin{aligned}\label{c6}
(1+(-1)^{p+q})R_{\bar p,\bar q}=&(-1)^p\sum_{j=0}^p (1+(-1)^j)\binom{p+q-j-1}{q-1}\zeta(j)R(\overline{p+q-j})\\
&+(-1)^p\sum_{k=0}^{q-1} (1-(-1)^k)\binom{p+q-k-2}{p-1}{R}({k+1}){R}(\overline{p+q-k-1})\\
&-(-1)^p(1-(-1)^q)\zeta(\bar p)R(\bar q),
\end{aligned}\\
&\begin{aligned}\label{c7}
(1+(-1)^{p+q})R_{p,\bar q}=&-(-1)^p\sum_{j=0}^p (1+(-1)^j)\binom{p+q-j-1}{q-1}\zeta(\bar j)R(\overline{p+q-j})\\
&-(-1)^p\sum_{k=0}^{q-1} (1+(-1)^k)\binom{p+q-k-2}{p-1}{R}(\overline{k+1}){R}({p+q-k-1})\\
&-(-1)^p(1-(-1)^q)\zeta(p)R(q),
\end{aligned}
\end{align}
where $\zeta(0)=\zeta(\bar 0):=-1/2,\ \zeta(1):=0$ and $R(1):=2\log(2)$.
\end{cor}

Note that the formula (\ref{c4}) can also be found in \cite{XZ2020}. From (\ref{c5})-(\ref{c7}), we can get the following cases.
\begin{exa} We have
\begin{align*}
&R_{{\bar 1},4}=-62\z(5)-\frac1{2}\pi R(\bar 4)-\frac1{4}\pi^3R(\bar 2)+15\log(2)\z(4),\\
&R_{{\bar 2},3}=93\z(5)+\frac7{2}\z(2)\z(3)+\frac3{2}\pi R(\bar 4)+\frac1{4}\pi^3R(\bar 2),\\
&R_{{\bar 1},{\bar 3}}=\frac3{2}R(\bar 4)-3\z(2)R(\bar 2)+\frac 1{4}\log(2)\pi^3,\\
&R_{{\bar 2},{\bar 2}}=-\frac3{2}R(\bar 4)+4\z(2)R(\bar 2),\\
&R_{1,\bar 3}=-\frac3{2}R(\bar 4)-\frac7{2}\pi\z(3)-\frac1{2}\log(2)\pi^3,\\
&R_{2,\bar 2}=\frac3{2}R(\bar 4)+\frac1{2}\z(2)R(\bar 2)+7\pi \z(3).
\end{align*}
\end{exa}

We see that in (\ref{c4}) and (\ref{c5}), if $p+q$ is odd, two modified forms of the identities hold, but without any (alternating) linear Euler $R$-sum occurring. Similarly, in (\ref{c6}) and (\ref{c7}), if $p+q$ is even, two modified forms of the identities hold, but without any (alternating) linear Euler $R$-sum occurring. These give back well-known nonlinear relations between (alternating) zeta values and $R$-values.  Next, in a same way, we establish a `duality' sum formula of Euler type sums.

\begin{thm}\label{thm3.3} For positive integers $m$ and $p$,
\begin{align}\label{c8}
&(-1)^m \sum_{i+j=m-1,\atop i,j\geq 0} \binom{p+i-1}{i}\binom{q+j-1}{j} \su \frac{M\B_n(p+i)}{(n-1/2)^{q+j}}a_n\nonumber\\
&+(-1)^p \sum_{i+j=p-1,\atop i,j\geq 0} \binom{m+i-1}{i}\binom{q+j-1}{j} \su \frac{M\a_{n-1}(m+i)}{(n-1/2)^{q+j}}b_{n-1}\nonumber\\
&=(-1)^{p+m-1}\binom{p+q+m-2}{q-1} \widehat{R}\ab(p+q+m-1)\nonumber\\
&\quad+(-1)^{p+m}\sum_{i+j=q-1,\atop i,j\geq 0} \binom{m+i-1}{i}\binom{p+j-1}{j}\widehat{R}\a(m+i)\widehat{R}\B(p+j).
\end{align}
\end{thm}
\pf Consider
\[f_2(s;A,B):=\frac{\Psi^{(m-1)}(-s;A)\Psi^{(p-1)}(-s;B)}{(m-1)!(p-1)!(s+1/2)^q},\]
which has poles of order $p+m$ at $s=-1/2$ and $n$ ($n$ is any non-negative integer). With the help of Lemma \ref{lem2.1}, the residues are easily calculated to be
\begin{align*}
{\rm Res}[f_2(s;A,B),s=n]&=(-1)^{p+m-1}\binom{p+q+m-2}{q-1}\frac{a_nb_n}{(n+1/2)^{p+q+m-1}}\\
&\quad-(-1)^m \sum_{i+j=m-1,\atop i,j\geq 0} \binom{p+i-1}{i}\binom{q+j-1}{j} \frac{M\B_n(p+i)}{(n+1/2)^{q+j}}a_n\\
&\quad-(-1)^p \sum_{i+j=p-1,\atop i,j\geq 0} \binom{m+i-1}{i}\binom{q+j-1}{j} \frac{M\a_n(m+i)}{(n+1/2)^{q+j}}b_n.
\end{align*}
Clearly, $f_2(s;A,B)$ also has a pole of order $p+q+m$ at $s=0$. Using (\ref{c2}) we find that
\begin{align*}
&{\rm Res}[f_2(s;A,B),s=0]=(-1)^{p+m}\sum_{i+j=q-1,\atop i,j\geq 0} \binom{m+i-1}{i}\binom{p+j-1}{j}\widehat{R}\a(m+i)\widehat{R}\B(p+j).
\end{align*}
Summing these two contributions, we thus immediately deduce (\ref{c8}) to complete the proof.\hfill$\square$

Hence, setting $A,B\in\{A_1,A_2\}$ in Theorem \ref{thm3.3} yields many linear relations between (alternating) linear Euler $R$-sums and polynomials in (alternating) single $R$-values. For example,
\begin{align*}
3R_{2,4}+2R_{3,3}=112\zeta^2(3)-\frac{\pi^6}{6}.
\end{align*}

Now, we evaluate the quadratic Euler type sums in the same manner as in the above.  Let $C:=\{c_k\}\ (-\infty < k < \infty)$ be a sequence of complex numbers with ${c_k} = o\left( {{k^\lambda }} \right)\ (\lambda < 1)$ if $k\rightarrow \pm \infty$, and let
\[\widehat{R}\abc(q):=\sk \frac{a_{k-1}b_{k-1}c_{k-1}}{(k-1/2)^q}.\]

\begin{thm}\label{thm3.4} Let $m,p$ and $q>1$ be positive integers with $A,B$ and $C$ defined above, we have
\begin{align}\label{c9}
&(-1)^{p+q+m} \su \frac{{\bar M}\B_n(m){\bar M}\C_n(p)}{(n-1/2)^q} a_n+\su \frac{M\B_{n-1}(m)M\C_{n-1}(p)}{(n-1/2)^q}a_{n-1}\nonumber\\
&+(-1)^{p+m}\binom{p+q+m-1}{q-1} \widehat{R}\abc(p+q+m)\nonumber\\
&+(-1)^m\sum_{j=1}^{m+1} \binom{j+p-2}{p-1}\binom{m+q-j}{q-1}\su \frac{M\C_{n-1}(j+p-1)}{(n-1/2)^{m+q-j+1}}a_{n-1}b_{n-1}\nonumber\\
&+(-1)^p\sum_{j=1}^{p+1} \binom{j+m-2}{m-1}\binom{p+q-j}{q-1}\su \frac{M\B_{n-1}(j+m-1)}{(n-1/2)^{p+q-j+1}}a_{n-1}c_{n-1}\nonumber\\
&-(-1)^{p+m} \sum_{j=1}^{p+m} \binom{p+q+m-j-1}{q-1} \su \frac{T\a_{n-1}(j)}{(n-1/2)^{p+q+m-j}}b_{n-1}c_{n-1}\nonumber\\
&-(-1)^m \sum_{j_1+j_2\leq m+1,\atop j_1,j_2\geq 1} \binom{m+q-j_1-j_2}{q-1}\binom{j_2+p-2}{p-1} \su \frac{T\a_{n-1}(j_1)M\C_{n-1}(j_2+p-1)}{(n-1/2)^{m+q-j_1-j_2+1}}b_{n-1}\nonumber\\
&-(-1)^p \sum_{j_1+j_2\leq p+1,\atop j_1,j_2\geq 1} \binom{p+q-j_1-j_2}{q-1}\binom{j_2+m-2}{m-1} \su \frac{T\a_{n-1}(j_1)M\B_{n-1}(j_2+m-1)}{(n-1/2)^{p+q-j_1-j_2+1}}c_{n-1}\nonumber\\
&+{\rm Res}[f_3(s;A,B,C),s=-1/2]=0,
\end{align}
where
\begin{align}\label{c10}
&{\rm Res}[f_3(s;A,B,C),s=-1/2]\nonumber\\
&=(-1)^{m+p}\sum_{k_1+k_2+k_3=q-1,\atop k_1,k_2,k_3\geq 0} \binom{m+k_2-1}{k_2}\binom{p+k_3-1}{k_3}\left((-1)^{k_1}R\a(k_1+1)-\widehat{R}\a(k_1+1) \right)\nonumber\\
&\quad\quad\quad\quad\quad\quad\quad\quad\quad\quad\times\widehat{R}\B(m+k_2)\widehat{R}\C(p+k_3).
\end{align}
\end{thm}
\pf Consider
\[f_3(s;A,B,C):=\frac{\pi\cot(\pi s;A)\Psi^{(m-1)}(-s;B)\Psi^{(p-1)}(-s;C)}{(m-1)!(p-1)!(s+1/2)^q},\]
and use residue computations to deduce the desired evaluation.\hfill$\square$

\begin{cor} For positive integers $p,m$ and $q>1$, if $p+q+m$ is even, then the three (alternating) quadratic Euler $R$-sums
\[R_{mp,q}=\sum_{n=0}^\infty \frac{H_{n}^{(m)}H_{n}^{(p)}}{(n+1/2)^q},\quad R_{m{\bar p},q}=\sum_{n=0}^\infty \frac{H_{n}^{(m)}{\bar H}_{n}^{(p)}}{(n+1/2)^q},\quad R_{{\bar m}{\bar p},q}=\sum_{n=0}^\infty \frac{{\bar H}_{n}^{(m)}{\bar H}_{n}^{(p)}}{(n+1/2)^q}\]
are reducible to (alternating) linear $R$-sums. And if $p+q+m$ is odd, then the three (alternating) quadratic Euler $R$-sums
\[R_{mp,{\bar q}}=\sum_{n=0}^\infty \frac{H_{n}^{(m)}H_{n}^{(p)}}{(n+1/2)^q}(-1)^{n},\ R_{m{\bar p},\bar q}=\sum_{n=0}^\infty \frac{H_{n}^{(m)}{\bar H}_{n}^{(p)}}{(n+1/2)^q}(-1)^{n},\ R_{{\bar m}{\bar p},\bar q}=\sum_{n=0}^\infty \frac{{\bar H}_{n}^{(m)}{\bar H}_{n}^{(p)}}{(n+1/2)^q}(-1)^{n}\]
are reducible to (alternating) linear $R$-sums.
\end{cor}
\pf Letting $A,B,C\in\{A_1,A_2\}$ in Theorem \ref{thm3.4} yields the desired description. \hfill$\square$

Note that the evaluation of $R_{mp,q}$ was proved in \cite{XZ2020}.

Next, we give a general reduction description of Euler $R$-sums.

\begin{thm} For positive integers $p_1,\ldots,p_r$ and $q>1$, (alternating) Euler $R$-sum $R_{e_1e_2\cdots e_r,q}$ $(e_j\in\{p_j,{\bar p}_j\}$ reduces to a combination of sums of lower orders whenever the weight $p_1+p_2+\cdots+p_r+q$ and the order $r$ are of the same parity. (alternating) Euler $R$-sum $R_{e_1e_2\cdots e_r,{\bar q}}$ $(e_j\in\{p_j,{\bar p}_j\}$ reduces to a combination of sums of lower orders whenever the weight $p_1+p_2+\cdots+p_r+q$ and the order $r$ are of the different parity.
\end{thm}
\pf Considering
\begin{align*}
\frac{\cot(\pi s; A)\Psi^{(p_1-1)}(-s;A^{(1)})\Psi^{(p_2-1)}(-s;A^{(2)})\cdots \Psi^{(p_r-1)}(-s;A^{(r)})}{(s+1/2)^q(p_1-1)!(p_2-1)!\cdots (p_r-1)!},\quad (A,A^{(l)}\in\{A_1,A_2\})
\end{align*}
and using residue computation, we obtain
\begin{align*}
&(-1)^{p_1+\cdots+p_r+q}\su \frac{\prod_{i=1}^r {\bar M}^{(A^{(i)})}_n(p_{i})}{(n-1/2)^q}a_n+\su \frac{\prod_{i=1}^r { M}^{(A^{(i)})}_{n-1}(p_{i})}{(n-1/2)^q}a_{n-1}\\
&\quad+(\text{sums of degree} \leq r)=0.
\end{align*}
Thus, with direct calculations, we may easily deduce the desired description.\hfill$\square$

\section{Evaluations of Multiple $R$-Values for lower depths}\label{sec4}
In this section, we establish some explicit evaluations of MRVs.

It is well-known (\cite{BBBL1997}) that the generating function of the `height one' multiple zeta values is given in terms of the gamma function:
\begin{align*}
1-\sum_{m,n=1}^\infty \z(m+1,\underbrace{1,\ldots,1}_{n-1})x^my^n=\frac{\Gamma(1-x)\Gamma(1-y)}{\Gamma(1-x-y)}.
\end{align*}
We can give the following $R$-version of this formula.
\begin{thm} We have the generating series identity
\begin{align}
&\sum_{m,n=1}^\infty {R}(m+1,\underbrace{1,\ldots,1}_{n-1})x^my^n=\frac{\Gamma(1/2)\Gamma(1-y)}{\Gamma(1/2-y)}-\frac{\Gamma(1/2-x)\Gamma(1-y)}{\Gamma(1/2-x-y)},
\end{align}
where $|x|<1/2,|y|<1/2$ and $|x+y|<1/2$.
\end{thm}
\pf According to the definition of MRVs, it is easy to see that
\begin{align*}
{R}(m+1,\underbrace{1,\ldots,1}_{n-1})=\frac{(-1)^{m+n-1}}{m!(n-1)!}\int_{0}^1 \frac{\log^m(t)\log^{n-1}(1-t)}{1-t}t^{-1/2}dt.
\end{align*}
Multiplying both sides by $x^my^n$ and summing both sides over $m,n$ gives
\begin{align*}
&\sum_{m,n=1}^\infty {R}(m+1,\underbrace{1,\ldots,1}_{n-1})x^my^n\\
&=\int_{0}^1 \sum_{m=1}^\infty \frac{1}{m!}(-x\log(t))^m\sum_{n=1}^\infty \frac{1}{(n-1)!}(-y\log(1-t))^{n-1} y\frac{t^{-1/2}}{1-t}dt\\
&=y\int_{0}^1 (t^{-x}-1)(1-t)^{-y-1}t^{-1/2}dt\\
&=yB(1/2-x,-y)-yB(1/2,-y).
\end{align*}
Then, we use the relation between Eulerian Beta function and Gamma function
\begin{align*}
B\left( {\alpha,\beta} \right) := \int\limits_0^1 {{x^{\alpha - 1}}{{\left( {1 - x} \right)}^{\beta - 1}}dx}  = \frac{{\Gamma \left( \alpha \right)\Gamma \left( \beta \right)}}{{\Gamma \left( {\alpha + \beta} \right)}}\quad ( {\mathop{\Re}\nolimits} \left( \alpha \right) > 0,{\mathop{\Re}\nolimits} \left( \beta \right) > 0)
\end{align*}
to complete the proof.

\begin{cor} For $|x|<1/2,|y|<1/2$ and $|x+y|<1/2$, it holds
\begin{align}\label{d2}
&\sum_{m,n=1}^\infty {R}(m+1,\underbrace{1,\ldots,1}_{n-1})x^my^n\nonumber\\
&=\exp\left(-2\log(2)y+2\sum_{n=2}^\infty \frac{(1-2^{n-1})y^n}{n}\zeta(n)\right)\nonumber\\
&\quad-\exp\left(-2\log(2)y+\sum_{n=2}^\infty \frac{(2^n-1)x^n+y^n-(2^n-1)(x+y)^n}{n}\zeta(n)\right).
\end{align}
\end{cor}
\pf Using the expression of Gamma function
\begin{align*}
\Gamma(1-x)=\exp\left(\gamma x+\sum_{n=2}^\infty \frac{\z(n)}{n}x^n\right)\quad(|x|<1).
\end{align*}
and the duplication formula $\Gamma(x)\Gamma(x+1/2)=2^{1-2x}\sqrt{\pi}\Gamma(2x)$, one obtain the expansion of $\log\Gamma(1/2-x)$:
\begin{align*}
\log\Gamma(1/2-x)=\frac{\log\pi}{2}+\gamma x+2x\log(2)+\sum_{n=2}^\infty \frac{(2^n-1)\z(n)}{n}x^n.
\end{align*}
Hence, an elementary calculation yields
\begin{align*}
\frac{\Gamma(1/2-x)\Gamma(1-y)}{\Gamma(1/2-x-y)}=\exp\left(-2\log(2)y+\sum_{n=2}^\infty \frac{(2^n-1)x^n+y^n-(2^n-1)(x+y)^n}{n}\zeta(n)\right).
\end{align*}
Thus, formula (\ref{d2}) holds. \hfill$\square$

Therefore, from (\ref{d2}), we have
\[{R}(m+1,\underbrace{1,\ldots,1}_{n-1})\in\mathbb{Q}[\log(2),\z(2),\z(3),\z(4),\ldots].\]
For example,
\begin{align*}
&R(2,1)=7\z(3)-6\z(2)\log(2),\\
&R(3,1)=\frac{45}{4}\z(4)-14\log(2)\z(3),\\
&R(2,1,1)=\frac{15}{2}\z(4)+6\log^2(2)\z(2)-14\log(2)\z(3),\\
&R(3,1,1)=62\z(5)-\frac{45}{2}\log(2)\z(4)-28\z(2)\z(3)+14\log^2(2)\z(3),\\
&R(2,1,1,1)=31\z(5)-15\log(2)\z(4)-13\z(2)\z(3)+14\log^2(2)\z(3)-4\log^3(2)\z(2).
\end{align*}

According to the definitions of MRVs and Euler $R$-sums, we have the relations
\begin{align*}
R(q,p)=R_{p,q},\ R(q,{\bar p})=-R_{{\bar p},q},\ R(\bar q,\bar p)=R_{\bar p,\bar q}\ \text{and}\ R(\bar q,p)=-R_{p,\bar q}.
\end{align*}
Hence, from identities (\ref{c4})-(\ref{c7}), we know that if $p+q$ is odd, then the (alternating) linear $R$-sums $R(q,p)$ and $R(q,{\bar p})$ are reducible to (alternating) zeta values and $R$-values, and if $p+q$ is even, then the (alternating) linear $R$-sums $R(\bar q,\bar p)$ and $R(\bar q,p)$ are reducible to (alternating) zeta values and $R$-values.

Further, by the methods of \cite{XW2018,H2016}, we may easily deduce the following relation
\begin{align*}
R_{i_1i_2\cdots i_m,q}
    =\sum_{\xi\in\mathcal{C}_m}\sum_{\sigma\in\mathcal{S}_m}
        \frac{R(q,J_1(I_{\sigma}^{(m)}),J_2(I_{\sigma}^{(m)}),\ldots,J_p(I_{\sigma}^{(m)}))}{\xi_1!\xi_2!\cdots\xi_p!}
       ,
\end{align*}
where $\xi:=(\xi_1,\xi_2,\ldots,\xi_p)\in\mathcal{C}_m$ ($\mathcal{C}_m$ is a set of all compositions of $m$) and a permutation $\sigma\in\mathcal{S}_m$ ($\mathcal{S}_m$ is a symmetric group of all permutations on $m$ symbols), $I_{\sigma}^{(m)}=(i_{\sigma(1)},\ldots,i_{\sigma(m)})$, and
\[
J_c(I_{\sigma}^{(m)})=i_{\sigma(\xi_1+\cdots+\xi_{c-1}+1)}+\cdots+i_{\sigma(\xi_1+\cdots+\xi_c)}\,,
    \quad\text{for } c=1,2,\ldots,p\,.
\]

Next, we give an explicit evaluations for multiple $R$-values with depth three. First,
from definitions,
\begin{align}\label{d3}
\sum_{n=1}^\infty \frac{h_n^{(k_1)}H_{n-1}^{(k_3)}}{n^{k_2}}=R(k_1)\z(k_2,k_3)-R(k_1,k_2,k_3)\quad (k_1,k_2>1),
\end{align}
where $h_n^{(p)}$ stands for the odd harmonic number, which is defined for $n\in\N_0$ and $p\in\N$ by
\begin{align*}
h_n^{(p)}:=\sum\limits_{k=1}^n \frac{1}{(k-1/2)^p},\quad h_0^{(p)}:=0,\quad h_n:=h_n^{(1)}.
\end{align*}
If $p>1$, the generalized odd harmonic number $h^{(p)}_n$ converges to the $R$-value:
\begin{align*}
\lim\limits_{n\rightarrow \infty}h^{(p)}_n=R(p).
\end{align*}

Now, we use the method of contour integration to prove an explicit reduction formula for the Euler type sum on the left-hand sides of (\ref{d3}).
By direct calculations, we can get the following lemma.
\begin{lem}\label{lem4.3} Let $p\geq 1$ and $n$ be nonnegative integers, if $|s-n+1/2|<1$, then
\begin{align}
\frac{\psi^{(p-1)}(-s)}{(p-1)!}=(-1)^p\sk \binom{k+p-2}{p-1}\left(R(k+p-1)+(-1)^{k+p-1}h_n^{(k+p-1)} \right)(s+1/2-n)^{k-1},
\end{align}
and if $|s+n-1/2|$ and $n>0$, then
\begin{align}
\frac{\psi^{(p-1)}(-s)}{(p-1)!}=(-1)^p\sk \binom{k+p-2}{p-1}\left(R(k+p-1)-h_{n-1}^{(k+p-1)} \right)(s-1/2+n)^{k-1},
\end{align}
where $R(1):=2\log(2)$ and $\z(1):=0$ wherever it occurs. When $p=1$, in which cases one should
replace $\psi^{(0)}(-s)$ by $\psi(-s)+\gamma$.
\end{lem}

Applying Lemma \ref{lem4.3}, we can prove the following theorem.

\begin{thm} For positive integers $p,m$ and $q>1$,
\begin{align}\label{d6}
&(1+(-1)^{p+q+m})\su\frac{H_{n-1}^{(m)}h_n^{(p)}}{n^q}\nonumber\\
&=-(-1)^{m+p}(1+(-1)^q)\z(m)\z(q)R(p)-(-1)^m(1-(-1)^{p+q})\z(m)\widetilde{T}(q,p)\nonumber\\
&\quad-(-1)^p(1-(-1)^{q+m})R(p)\z(q,m)-(-1)^pR(p)\z(q+m)-\widetilde{T}(m+q,p)\nonumber\\
&\quad+(-1)^{m+p}\sum_{k=1}^{m+1} (-1)^k \binom{k+p-2}{p-1} \binom{m+q-k}{q-1}\left\{R(k+p-1)\z(m+q-k+1)\atop+(-1)^{k+p-1}\widetilde{T}(m+q-k+1,k+p-1) \right\}\nonumber\\
&\quad-2(-1)^{m+p}\sum_{2k_1+k_2\leq m+1,\atop k_1,k_2\geq 1} (-1)^{2k_1+k_2}\binom{k_2+p-2}{p-1}\binom{m+q-2k_1-k_2}{q-1}\z(2k_1)\nonumber\\
&\quad\quad\quad\quad\quad\quad\quad\quad\quad\quad\times\left\{R(k_2+p-1)\z(m+q-2k_1-k_2+1)\atop+(-1)^{k_2+p-1}\widetilde{T}(m+q-2k_1-k_2+1,k_2+p-1)\right\}\nonumber\\
&\quad+(-1)^m\sum_{k_1+k_2+k_3=p-1,\atop k_1,k_2,k_3\geq 0}(1-(-1)^{k_1})(-1)^{k_3}\binom{k_2+m-1}{k_2}\binom{k_3+q-1}{k_3}R(k_1+1)\nonumber\\
&\quad\quad\quad\quad\quad\quad\quad\quad\quad\quad\times\left\{R(m+k_2)R(k_3+q)+(-1)^{m+k_2}\t(k_3+q,m+k_2)\right\}\nonumber\\
&\quad-{\rm Res}[F(s),s=0],
\end{align}
where $R(1):=2\log(2)$ and $\z(1):=0$ wherever it occurs, and
\begin{align}\label{d7}
&{\rm Res}[F(s),s=0]\nonumber\\
&=(-1)^p\binom{p+q+m-1}{p-1} R(p+q+m)\nonumber\\ &\quad-2(-1)^p \sum_{k=1}^{[(m+q)/2]}\binom{p+q+m-2k-1}{p-1}\z(2k)R(p+q+m-2k)\nonumber\\
&\quad+(-1)^{m+p}\sum_{k=1}^{q+1} \binom{k+m-2}{m-1}\binom{p+q-k}{p-1} \z(k+m-1)R(p+q-k+1)\nonumber\\
&\quad-2(-1)^{m+p}\sum_{2k_1+k_2\leq q+1,\atop k_1,k_2\geq 1} \binom{k_2+m-2}{m-1}\binom{p+q-2k_1-k_2}{p-1}\nonumber\\&\quad\quad\quad\quad\quad\quad\quad\quad\quad\quad\times\z(2k_1)\z(k_2+m-1)R(p+q-2k_1-k_2+1).
\end{align}
\end{thm}
\pf  Consider
\begin{align*}
F(s):=\frac{\pi\cot(\pi s)\psi^{(m-1)}(-s)\psi^{(p-1)}(-s)}{s^q(m-1)!(p-1)!}
\end{align*}
which has poles at $s=0,\pm n$ and $n-1/2$ ($n$ is any positive integer). At a negative integer $-n$ the pole is simple and
the residue is
\begin{align*}
{\rm Res}[F(s),s=-n]=(-1)^{p+q+m} \frac{\left(\z(m)-H^{(m)}_{n-1}\right)\left(R(p)-h_n^{(p)}\right)}{n^q}.
\end{align*}
At a positive integer $n$, the pole has order $m+1$ and the residue is
\begin{align*}
{\rm Res}[F(s),s=n]&=(-1)^{m+p}\sum_{k=1}^{m+1} (-1)^{k-1}\binom{k+p-2}{p-1}\binom{m+q-k}{q-1}\\&\quad\quad\quad\quad\quad\quad\quad\times \frac{R(k+p-1)+(-1)^{k+p-1}h_n^{(k+p-1)}}{n^{m+q-k+1}}\\
&\quad+2(-1)^{m+p}\sum_{2k_1+k_2\leq m+1,\atop k_1,k_2\geq 1} (-1)^{2k_1+k_2}\binom{k_2+p-2}{p-1}\binom{m+q-2k_1-k_2}{q-1}\\&\quad\quad\quad\quad\quad\quad\quad\quad\quad\times\z(2k_1) \frac{R(k_2+p-1)+(-1)^{k_2+p-1}h_n^{(k_2+p-1)}}{n^{m+q-2k_1-k_2+1}}\\
&\quad+(-1)^{m+p}\frac{\left(\z(m)+(-1)^mH^{(m)}_n\right)\left(R(p)+(-1)^ph_n^{(p)} \right)}{n^q}.
\end{align*}
At a rational number $n-1/2$  the pole has order $p$ and the residue is
\begin{align*}
{\rm Res}[F(s),s=n-1/2]&=(-1)^m\sum_{k_1+k_2+k_3=p-1,\atop k_1,k_2,k_3\geq 0} ((-1)^{k_1}-1)(-1)^{k_3} \binom{k_2+m-1}{k_2}\binom{k_3+q-1}{k_3}\\&\quad\quad\quad\quad\quad\quad\quad\times R(k_1+1)\frac{R(m+k_2)+(-1)^{m+k_2}h_n^{(m+k_2)}}{(n-1/2)^{k_3+q}}.
\end{align*}
Finally, we can compute that the residue of the pole of order at $0$ is (\ref{d7}). Summing these four contributions yields the statement of the theorem. \hfill$\square$

Hence, from (\ref{d3}) and (\ref{d6}), we obtain the following description.
\begin{cor} For positive integers $k_1,k_2,k_3$ with $k_1>1$, the triple $R$-values $R(k_1,k_2,k_3)$ can be expressed in terms of combinations of double $R$-values, double $T$-values, double $t$-values, double zeta values and single zeta values.
\end{cor}

Since, the multiple $R$-values, multiple $T$-values and multiple $t$-values can be expressed as a $\mathbb{Q}$-linear combination of alternating multiple zeta values with same depth and weight, so we have
\begin{align*}
R(k_1,k_2,k_3)\in\mathbb{Q}[\text{(alternating) double zeta values and single zeta values}.]
\end{align*}

We end this paper by proposing the following conjecture.

\begin{con}
Let ${\bf k}=(k_1,\ldots,k_r)$ be an admissible index and assume its depth $r$ and weight $k_1+\cdots+k_r$ are of different parity. Then $R({\bf k})$ can be expressed as a $\mathbb{Q}$-linear combination of (alternating) multiple zeta values of lower depths and products of (alternating) multiple zeta values with sum of depths not exceeding $r$?
\end{con}
{\bf Acknowledgments.}  The author expresses his deep gratitude to Professors Masanobu Kaneko and Jianqiang Zhao for valuable discussions and comments.

 \end{document}